\documentclass[11pt]{amsart}
\usepackage[margin=1.2in]{geometry}

\usepackage{hyperref}

\usepackage[parfill]{parskip}

\usepackage{latexsym, amssymb, bbm}
\usepackage{amsthm}
\usepackage{xypic}
\usepackage{xcolor}
\usepackage{graphicx}
\usepackage{amsmath}
\usepackage{enumerate}

\newtheorem{lem}{Lemma}[section]
\newtheorem{thm}[lem]{Theorem}
\newtheorem{prop}[lem]{Proposition}
\newtheorem{cor}[lem]{Corollary}

\newtheorem{remark}[lem]{Remark}

\theoremstyle{definition}
\newtheorem{defn}[lem]{Definition}

\newcommand{\wt}[1]{\widetilde{#1}}

\newcommand{\bdf}{\begin{defn}}
\newcommand{\edf}{\end{defn}}

\newcommand{\bthm}{\begin{thm}}
\newcommand{\ethm}{\end{thm}}

\newcommand{\blem}{\begin{lem}}
\newcommand{\elem}{\end{lem}}

\newcommand{\bcor}{\begin{cor}}
\newcommand{\ecor}{\end{cor}}

\newcommand{\bprop}{\begin{prop}}
\newcommand{\eprop}{\end{prop}}

\newcommand{\brmk}{\begin{remark}}
\newcommand{\ermk}{\end{remark}}

\newcommand{\bpf}{\begin{proof}}
\newcommand{\epf}{\end{proof}}

\newcommand{\beq}{\begin{equation}}
\newcommand{\eeq}{\end{equation}}

\newcommand{\bit}{\begin{itemize}}
\newcommand{\eit}{\end{itemize}}

\newcommand{\ben}{\begin{enumerate}}
\newcommand{\een}{\end{enumerate}}

\numberwithin{equation}{section}

\def\C{\mathbb{C}}

\def\R{\mathbb{R}}
\def\Z{\mathbb{Z}}

\def\CP{\mathbb{CP}}
\def\RP{\mathbb{RP}}

\def\cB{\mathcal{B}}
\def\CC{\mathcal{C}}

\def\cE{\mathcal{E}}
\def\FF{\mathcal{F}}

\def\JJ{\mathcal{J}}
\def\cJ{\mathcal{J}}

\def\NN{\mathcal{N}}
\def\cO{\mathcal{O}}

\def\cS{\mathcal{S}}

\def\w{\omega}

\def\xkm2{\overline{X}_{k-2}}

\begin{document}
\title{Exact Lagrangians in $A_n$-surface singularities}



\author[Wu]{Weiwei Wu}

\thanks{Address: Department of Mathematics, Michgan State University, USA}
\thanks{Email: wwwu@math.msu.edu}
\thanks{The author is supported by NSF Focused Research Grants DMS-0244663.}

\date{\today}


\maketitle

\begin{abstract}
In this paper we classify Lagrangian spheres in $A_n$-surface singularities up to Hamiltonian
isotopy.  Combining with a result of A. Ritter \cite{Ri10}, this yields a complete classification of
exact Lagrangians in $A_n$-surface singularities.  Our main new tool is the application of
the a technique which we call \textit{ball-swappings} and its relative version.
\end{abstract}

\bigskip
\noindent
{\bf MSC classes:}  53Dxx, 53D35, 53D12

\bigskip
\noindent
{\bf Keywords:} symplectic ball packing, Lagrangian isotopy, symplectomorphism groups

\section{Introduction}

One of the classical problems in sympletic topology is to understand the classification of exact Lagrangians
in a given exact symplectic manifold.  As appealing as it is, the problem is in general very difficult, even
in its simplest form.  In particular, Arnold's nearby Lagrangian conjecture asserts that any exact
Lagrangian in $T^*M$ is Hamiltonian isotopic to the zero section, which is still quite open.  This line
of questions have attracted many efforts involving a long list of authors, among which we only mention a few
very recent advances \cite{Ab12, AI12, Se12}.

Note that the proofs of above-mentioned works involved advanced techniques from Floer theory, and usually
covers symplectic manifolds of arbitrary dimensions.  While these deep methods are amazingly
powerful in determining the shape of exact Lagrangians coupling with homotopy methods, a general construction
 of Hamiltonian isotopies is still missing in dimension higher than $4$.
 In contrast, one has more geometric tools available
in dimension $4$, thus equivalence up to Hamiltonian classes is more approachable.  See \cite{Hi12} for an example of
 a very nice application of foliation techniques to this problem, and \cite{LW12, HPW12} for another approach which is
 closer in idea to what we will use here.

In the present paper, we investigate the classification of Lagrangian spheres in an $A_n$-surface singularity.
By definition, an $A_n$-surface singularity is symplectically identified with the subvariety
$$\{(x,y,z): x^2+y^2+z^{n+1}=1\}\subset (\C^3,\w_{std}),$$

endowed with the restricted K\"ahler form.  Throughout the paper we will denote $W$ as the
$A_{n-1}$-surface singularity, which is the main object that we will investigate.  It is
by now well-known that $W$ is identified symplectically with the plumbing of $n-1$ copies $T^*S^2$.  We will call the zero sections of these
plumbed copies \textit{standard spheres}.  The following is our main result:

\bthm\label{t:mainL} Lagrangian spheres in $W$ are unique up to Hamiltonian isotopy and Lagrangian Dehn
twists along the standard spheres. \ethm

A fantastic result showed by A. Ritter in \cite{Ri10} says that, embedded exact Lagrangians in $W$ are all Lagrangian spheres.
We therefore obtain the following corollary, which completely classifies exact Lagrangians in $A_n$ surface singularities
up to Hamiltonian isotopy:

\bcor\label{c:CEL} Exact Lagrangians in $A_n$-surface singularities are isotopic to the zero section of a
plumbed copy of $T^*S^2$, up to a composition of Lagrangian Dehn twists along the standard spheres.
\ecor

Such kind of classification seems desirable but rare in the literature, especially when there exist
smoothly isotopic but not Hamiltonian isotopic Lagrangians \cite{Se99}.
Notice that various forms of partial results have been obtained previously.  In particular, R. Hind in \cite{Hi12} proves
Theorem \ref{t:mainL} for the case of $A_1$ and $A_2$.  It was also known that the result is true up to equivalence
of objects in the Fukaya category, using the deep computations in algebraic geometry as well as the
mirror symmetry of $A_n$-surface singularities \cite{IU05,IUU10}.  This Floer-theoretic version already found
interesting applications \cite{LM12, Se12}.  In principle, Theorem \ref{t:mainL} along with computation of
\cite{KS02} on the symplectic side should recover corresponding results in \cite{IU05,IUU10} on the mirror side.

\textit{En route}, we also prove the following result on the compactly supported symplectomorphism group of $W$:

\bthm\label{t:mainS} Any compactly supported symplectomorphism is Hamiltonian isotopic to a composition of
Dehn twists along the standard spheres.  In particular, $\pi_0(Symp_c(W))=Br_n$. \ethm

This is a refinement of the results due to J. Evans \cite[Theorem 4]{Ev11}, which asserts that $\pi_0(Symp_c(W))$
injects into $Br_n$, and Khovanov-Seidel \cite[Corollary 1.4]{KS02}, which proves in any dimension of $A_n$-singularities,
 there is an injection from the other direction.  Our result shows that these two injections are in fact both isomorphisms
 in dimension 4.

The paper is structured as follows.  In Section \ref{s:setup} we set up the notation and describe two different but closely
related models for $W$ and its compactification.  We then reduce the main theorems to problems in its compactification.
Section \ref{s:MCGW} contains the proof of Proposition \ref{p:MainS}, which implies Theorem \ref{t:mainS}, and section \ref{s:LagU}
contains the proof of Proposition \ref{p:MainL}, which implies Theorem \ref{t:mainL}.  We conclude the article with some
discussions on the ball-swapping symplectomorphisms, which is the main technique involved in this article.

\subsubsection*{Acknowledgements}
The author is deeply indebted to Richard Hind, who inspired the idea in this paper during his stay in Michigan State University.
I warmly thank Dusa McDuff for explaining many details of her recent preprint \cite{Mc12}.  Many ideas involved in this work
were originally due to Jonny Evans in his excellent series of papers \cite{Ev10, Ev11}, and the ball-swapping construction used
here stems out from the beautiful ideas exhibited in Seidel's lecture notes \cite{Se08}.
I would also like to thank the Geometry/Topology group of Michigan State University
 for providing me a friendly and inspiring working environment,
as well as supports for my visitors.

\section{Two models of $W$}\label{s:setup}

In this section we recall two models of $W$ and its compactifications due to J. Evans and I. Smith.  Along the way we introduce
various notations that we will use throughout the paper.

We first go through the compactification of the $A_n$-singularities as a complex affine variety following \cite{Ev11}.

Let $M$ be the blow-up of $\CP^2$ at $\{p_i=[\xi_i,0,1]\}_{i=1}^n$, which is identified as
$$\{([x,y,z], [a_i,b_i]): [x,y,z]\in \CP^2, [a_i,b_i]\in \CP^1, a_k y=b_k(x-\xi_k z), i=1,\dots,n\}.$$

Here $\xi_i$ denotes the $i$-th $n$-unit root. We therefore obtain a pencil structure of $M$, which is the blow-up of the pencil of lines passing through $[0,1,0]$ in $\CP^2$.  We denote:

\begin{itemize}
  \item the pencil as $\{P_t\}_{t=[x,z]\in\CP^1}$ ,
  \item $C_i$ is the exceptional curves of the blow-up at $[\xi_i,0,1]$, $i=1,\dots, n$,
  \item $C_{n+1}=\{[x,y,0], [x,y],\dots, [x,y]\}$,
  \item $C_{n+2}=\{[x,0,z], [0,1],\dots, [0,1]\}$.
\end{itemize}

Here $C_{n+1}$ plays the role of a generic fiber of the pencil which is a line, and $C_{n+2}$ is the proper transform of the line
passing through $\{p_i\}_{i=1}^n$, which is also a section of the pencil.  Endow a K\"ahler form $\w$ to $M$ so that
$\int_{C_{n+1}}\w=1$, $\int_{C_i}\w=r:=1/N$ for all $1\leq i\leq n$.  Here we require $N\gg n$ to be a large integer.  This
is always possible by the construction of a symplectic blow-up \cite{MS98}.

Let $U=M\backslash(C_{n+1}\cup C_{n+2})$.  Then $U$ has a Lefschetz fibration structure induced from $M$.
Lemma 7.1 of \cite{Ev11} showed that $U$ is biholomorphic to the $A_n$-singularity $W$.
As a result of \cite[Lemma 2.1.6]{EvT}, $U$ has a symplectic completion symplectomorphic to $W$.
Therefore, one may obtain a symplectic embedding $\iota: U\hookrightarrow W$. Note that the Lefschetz fibration of
$U$ defined above coincides with that of $W$.  In particular they have the same number of
Lefschetz singularities and monodromies.  Thus one may assume $\iota$ preserves the Lefschetz fibration structure.
In particular, their Lefschetz thimbles, thus matching cycles coincide.  If we identify $C_{n+2}\backslash C_{n+1}$ with
the base of the Lefschetz fibration of $U$, then the
standard spheres are the matching cycles lying above the straight arcs connecting $p_i$ to $p_{i+1}$ in the
base.  We will use this interpretation of standard spheres throughout, but a more explicit symplectic description is
in order.


\begin{figure}[h]
\centering
\includegraphics[trim = 0mm 25cm 13cm 0mm, clip, width=3in]{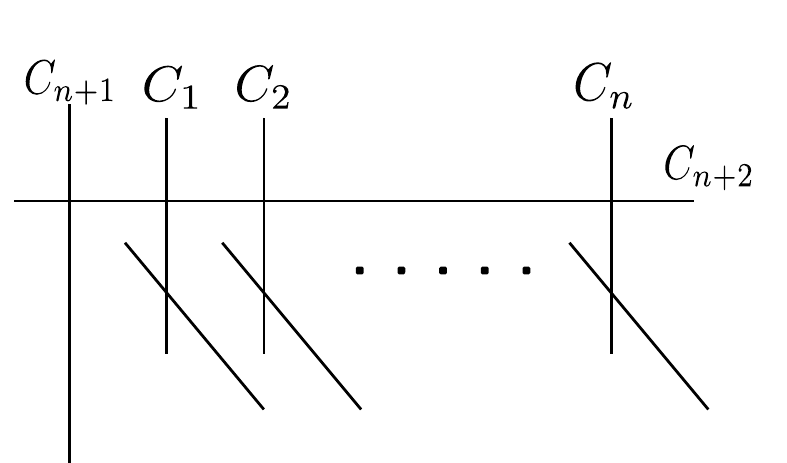}
\caption{Evans's Construction} \label{f:Evans}
\end{figure}

\begin{figure}[h]
\centering
\includegraphics[trim = 5mm 18cm 10mm 3cm, clip, width=6in]{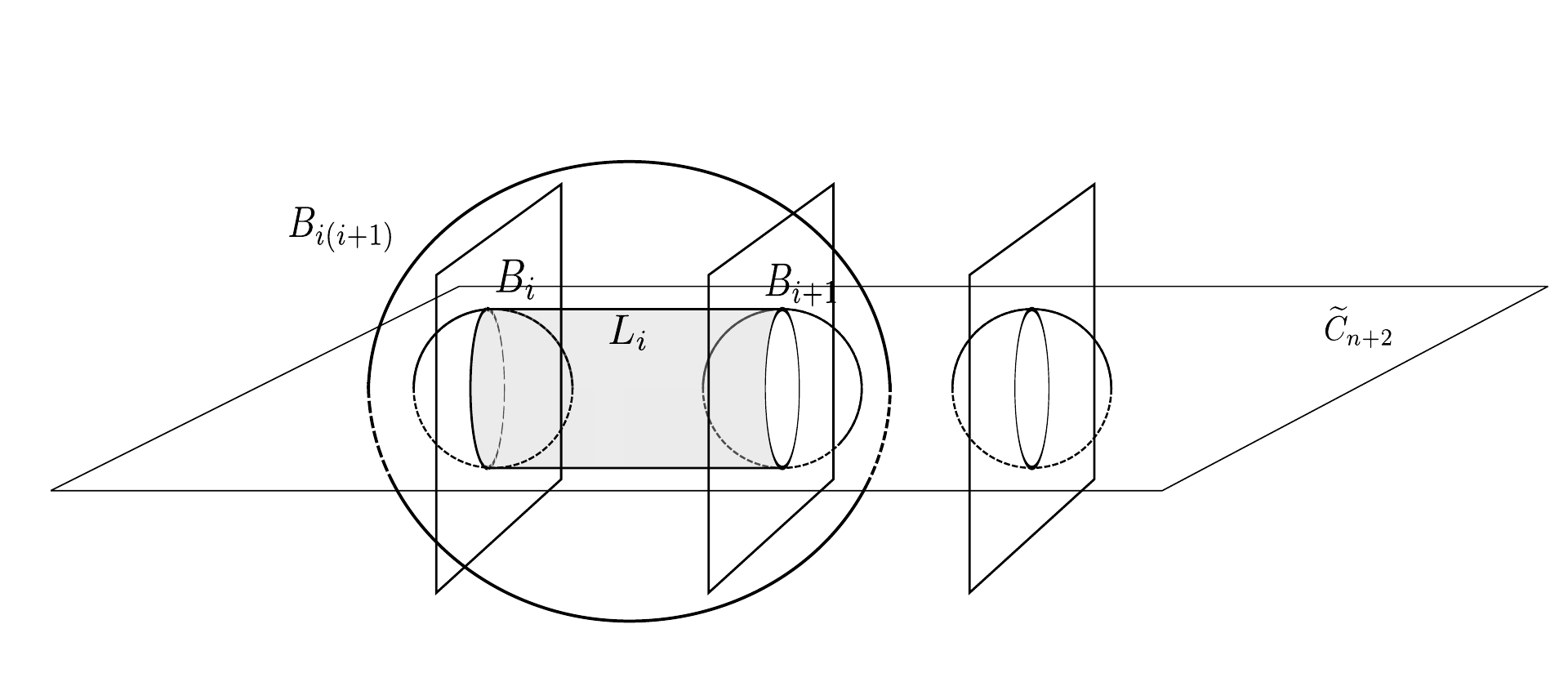}
\caption{Smith's Construction} \label{f:Smith}
\end{figure}


The above compactification model is explicit in complex coordinates,
but for our purpose of constructing symplectomorphisms
we need to recall the following alternative model slightly generalized
 from \cite[Example 4.25]{Sm12}.

For $z^0=(z^0_1,z^0_2)\in \C^2$ and real number $R$, let
$$B(z^0;R)=\{(z_1,z_2)\in \C^2: |z_1-z^0_1|^2+|z_2-z^0_2|^2<R^2\}.$$

Consider $B=B(0,\frac{1}{\sqrt{\pi}}),B_i=B(\frac{2i-n-1}{(n+1)\sqrt{\pi}},\frac{1}{\sqrt{N\pi}})$ for $1\leq i\leq n$.  The
circles $\{|z_2|=\frac{1}{\sqrt{N\pi}}\}$ over the arcs
$$\gamma_i=\{im(z_1)=0, Re(z_1)\in [\frac{2i-n-1}{(1+n)\sqrt{\pi}},\frac{2i-n+1}{(1+n)\sqrt{\pi}}]\}\subset \{z_2=0\},\hskip 2mm 1\leq i\leq n-1,$$

forms $(n-1)$ Lagrangian tubes. Upon blowing up all $B_i$, the boundary circles of the tubes are collapsed to points on the
exceptional spheres, and the Lagrangian tubes become
matching cycles between two singular fibers of the Lefschetz fibration described previously, thus giving the standard
spheres up to Hamiltonian isotopy.

The relation between Evans's and Smith's construction is clear from the symplectic interpretation of blow-ups.  In particular,
since the complement of a line in $\CP^2$ is identifies as a $4$-ball,
$U$ is exactly the complement of the proper transform of $\{z_2=0\}$ in the blow up of $B(0,\frac{1}{\sqrt{\pi}})$ along $B_i$.
Notice the following facts:

%



\blem[\cite{Ev11}, Proposition 2.1]\label{l:U=W} $Symp_c(U)$ is weakly homotopic to $Symp_c(W)$.

\elem

\blem\label{l:reduction} Suppose Lagrangian spheres are unique up to compactly supported symplectomorphisms in $U$, then the same
holds in $W$.
\elem

To see \ref{l:reduction}, notice that since $U$ has the symplectization identified with $W$, any
Lagrangian $L\subset W$ is isotopic to one in $U$ through the negative Liouville flow.
The two lemmata reduce our main theorems to the compactified case.
Concretely, they show that the following two propositions imply the main theorems:

\bprop\label{p:MainS} $\pi_0(Symp_c(U))$ is generated by the Dehn twists along $L_1,\dots, L_n$, where $L_i$ are matching
cycles of the Lefschetz fibration of $U$, for $1\leq i\leq n$.

\eprop

\bprop\label{p:MainL} Any pair of Lagrangian spheres $L,L'\subset U$ are symplectomorphic, i.e.,
there is a $\phi\in Symp_c(U)$ such that $\phi(L)=L'$.

\eprop



\section{The mapping class group of $W$}\label{s:MCGW}

In this section we introduce a technique of producing symplectomorphism alluded in \cite{LW12, BLW12},
which we call the \textit{ball-swapping}, and try to address its relation with the Dehn twists.

We start with a more general context.  Suppose $X$ is a symplectic manifold.  Given two symplectic ball embeddings:

$$\iota_{0,1}:\coprod_{i=1}^n B(r_i)\rightarrow X,$$

where $\iota_0$ is isotopic to $\iota_1$ through a Hamiltonian path $\{\iota_t\}$.  From the interpretation of blow-ups
in the symplectic category \cite{MP94}, the blow-ups can be represented as

$$X^{\#\iota_j}=(X\backslash\iota_j(\coprod_{i=1}^n B_i))/\sim,\text{ for }j=0,1.$$

Here the equivalence relation $\sim$ collapses the natural $S^1$-action on $\partial B_i=S^3$.  Now assume that
$K=\iota_0(\coprod B_i)=\iota_1(\coprod B_i)$ as sets, then $\iota_t$ defines a symplectic automorphism $\widetilde{\tau}_\iota$ of
$X\backslash K$, which descends to an automorphism $\tau_\iota$ of $X^{\#\iota}:=X^{\#\iota_0}=X^{\#\iota_1}$.  We call
$\tau_\iota$ a \textit{ball-swapping symplectomorphism} or \textit{ball-swapping} on $X^{\#\iota}$.  Notice it is not known to be true (nor false) that any two ball packings are Hamiltonian isotopic.  D. McDuff gave a comfirmative answer to this question for the symplectic $4$-manifolds with $b^+=1$ \cite{Mc98} which allows more freedom to create ball-swappings in the blow-ups of these manifolds, but the general case of the question is still widely open.  This construction is closely related to one in algebraic geometry, see discussions in Section $5$.

As usual, we denote $\tau_L$ the Lagrangian Dehn twist along $L$ for a Lagrangian sphere $L$ \cite{Se08}.
To compare the ball-swapping, the Dehn twists in dimension $4$ and the full mapping class group of $W$,
we first need a local refinement of the following result due to Evans:

\bthm[\cite{Ev11}, Theorem 1.4]\label{t:EvIn} The compactly supported symplectomorphism group $Symp_c(W)$ has weakly contractible connected components.  Moreover, $\pi_0(Symp_c(W))$ has an injective homomorphism into the pure braid group $Br_n$ with $n$-strands.

\ethm

We go over the main ingredients of Evans' proof briefly, which indeed proves the conclusion for $Symp_c(U)$.
Define a standard configuration $\{S_i\}_{i=1}^n$ in $M=\CP^2\#n\overline{\CP^2}$ as:

\ben[(i)]
   \item each $S_i$ is disjoint from $C_{n+1}$,
   \item $[S_i]=[C_i]$,
   \item there exist $J\in\mathcal{J}_\w$, the set of almost complex structures compatible with $\w$,
           for which all $S_i,C_{n+1},C_{n+2}$ are $J$-holomorphic.
   \item There is a neighborhood $\nu$ of $C_{n+2}$ such that $S_i\cap\nu=P_{t_i}\cap\nu$, for $t_i=S_i\cap C_{n+2}$.
\een

Proposition 7.4 of \cite{Ev11} showed that such configurations form a weakly contractible space $\CC_0$.
The proof needs to ensure no bubble occurs from $S_i$ to apply Pinsonnault's result \cite[Lemma 1.2]{Pi08}.
This partly motivates our choice of areas of $\omega(C_i)$.

Let $Conf(n)$ be the configuration space of $n$ points on a disk.  Proposition 7.2 of \cite{Ev11} shows one has the following homotopy fibration:
\begin{align*}
\Psi:\CC_0&\rightarrow Conf(n),\\
S=\bigcup_{i=1}^n S_i&\mapsto\{S_1\cap C_{n+2},\dots, S_n\cap C_{n+2}\}.
\end{align*}

Here the disk is identified with $C_{n+2}\backslash C_{n+1}$.  The fiber of this fibration is denoted as
$\FF$.  Explicitly, fix an unordered $n$-tuple of points
$[x_1,\dots,x_n]$, $x_i\in (D_2)$, $\FF$ is the space of standard configurations $\{S_i\}$ such that
$[S_i\cap C_{n+2}]_{i=1}^n=[x_i]_{i=1}^n$ as unordered $n$-tuples.
The associated long exact sequence thus gives the following isomorphism:

\beq\label{e:Br=Con} Br_n=\pi_1(Conf(n))\xrightarrow{\sim} \pi_0(\FF) \eeq

Moreover, \cite[Theorem 7.6]{Ev11} shows that:

\beq\label{e:ContF} \pi_i(\FF)=0, \text{ for all }i>1, \eeq

that is, the connected components of $\FF$ are weakly contractible.
The construction of the isomorphism \eqref{e:Br=Con} amounts to the following lemma:

\blem[\cite{Ev11}, Proposition 7.2]\label{l:EvE}
 Let $\alpha$ be a loop in $Conf(n)$.  One may construct a compactly supported Hamiltonian on $C_{n+2}\backslash C_{n+1}$ inducing
$\alpha$, then there is an extension of $\alpha$ to a Hamiltonian $f_\alpha$ of $M$ supported in a neighborhood of $C_{n+2}$, such that:

\bit
\item $f_\alpha$ preserves $C_{n+1}\cup C_{n+2}$ and fixes $C_{n+1}$ pointwisely,
\item $f_\alpha$ preserves the set $\CC_0$.

\eit

\elem

Take $[S_1,\dots, S_n]\in\FF$ and $\alpha$ a loop in $Conf(n)$.  Then \eqref{e:Br=Con} is given by
$[f_\alpha(S_1),\dots, f_\alpha(S_n)]$.  Note that there is by no means a canonical group structure on $\pi_0(\FF)$,
so one should understand \eqref{e:Br=Con} as the $Br_n$-action on $\pi_0(\FF)$ is free and transitive.
This explicit construction of the isomorphism \eqref{e:Br=Con} will be used later.

We now look closer to the special case when $n=2$ using Smith's model.
Let $B=B(0;1), B_+=B((\frac{1}{2},0);\frac{1}{3}), B_-=B((-\frac{1}{2},0);\frac{1}{3})$, and $\wt C=\{z_2=0\}\cap B$.
From the packing-blowup correspondence, $B\backslash (B_+\cup B_-\cup C)$ is exactly $U$ constructed previously for $n=2$, $N=\frac{1}{3}$.  We denote this open symplectic manifold by $U_2$ in case of confusions.  Therefore, by Lemma \ref{l:U=W} and \cite{Se98}
\beq\label{e:SympcU} Symp_c(U)=Symp_c(T^*S^2)=\Z,\eeq

and it is generated by the Lagrangian Dehn twist along the matching cycle described in Section \ref{s:setup}.
Now we consider a Hamiltonian isotopy $\rho$ swapping $B_-$ and $B_+$.  Explicitly, we choose a small number $\epsilon>0$, and
 a smooth function $f:\R\rightarrow\R$ such that $f(x)=0$ when $x\geq 1-\epsilon$ and $f(x)=1$ when $x\leq 1-2\epsilon$.  Then
 $\rho$ is defined by the Hamiltonian function $\pi f(|z_1|^2+|z_2|^2)\cdot|z_1|^2$.  $\rho$ thus defines a ball-swapping symplectomorphism
  $\tau_\rho\in Symp_c(B\#2\overline{\CP^2})$ with corresponding symplectic form.
  Denote $C$ as the proper transform of $\wt C$
  under this blow-up (which is $C_4$ in the notation of general cases), and $C_+$ ($C_-$ resp.) for the exceptional sphere
  by blowing up $\rho(B_+)$ ($\rho(B_-)$ resp.).  Then $\tau_\rho$ preserves $C$ and exchanges $C_+, C_-$ as sets.

Consider $l:=\tau_\rho|_{C}:(C, [x_+,x_-])\rightarrow(C, [x_+,x_-])$,
where $x_\pm=C_\pm\cap C$ and $[x_+,x_-]$ denotes the unorder pair of points they form on $\wt C$.
$l$ can be symplectically isotoped to $id$ on $C$ by a compactly supported Hamiltonian path $l_t^{-1}$.
We may then lift this isotopy to $l_t^{-1}\times id\in Symp_c(D_2(1)\times D_2(\epsilon))$ for some small $\epsilon>0$
and extend this lift to a symplectomorphism $\overline{l^{-1}}$.  From our way of constructing the symplectomorphisms,
it is clear that $\overline{l^{-1}}\circ\tau_\rho=id$ on a neighborhood of $C$.  Therefore,
$\overline{l^{-1}}\circ\tau_\rho\in Symp_c(U_2)$.

\begin{figure}[h]
\centering
\includegraphics[trim = 0mm 10cm 0mm 2cm, clip, width=6in]{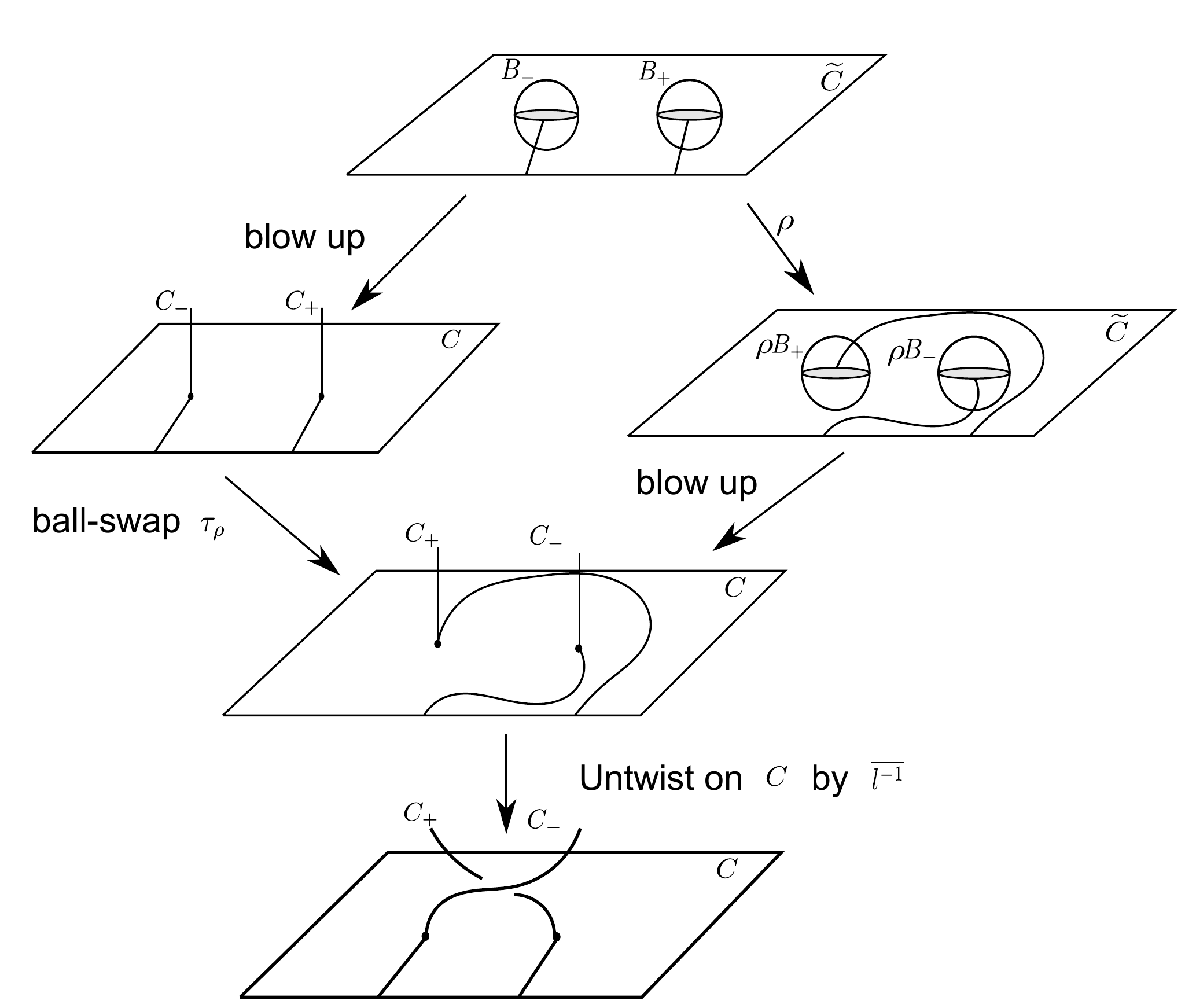}
\caption{Swapping two balls} \label{f:2ballswap}
\end{figure}

\blem\label{l:BDlocal}
$\overline{l^{-1}}\circ\tau_\rho$ is isotopic to a Lagrangian Dehn twist which generates $\pi_0(Symp_c(U_2))$.
\elem

\bpf

     Consider the action of $\overline{l^{-1}}\circ\tau_\rho$ on $([C_+, C_-])\in\FF$.  Since
     $\tau_\rho$ exchanges $C_+$ and $C_-$ as sets, the action of $\tau_\rho$ on the unordered pair $[C_+,C_-]$ is actually trivial.
     Also, notice that $l^{-1}([x_+,x_-])$ is the generator of $\pi_1(Conf(2))$, from \eqref{e:Br=Con} we see
     that  $\{(\overline{l^{-1}}\circ\tau_\rho)^k([C_+, C_-]): k\in\Z\}$ contains exactly one point in each component of
     $\FF$.  On the other hand, since $\tau_L$ is the generator of $\pi_0(Symp_c(U))$,
      $\overline{l^{-1}}\circ\tau_\rho$ is compactly isotopic to $\tau_L^m$ for some integer $m$.  Therefore,
       $\tau_L$ also acts transitively on $\pi_0(\FF)$.  Now \cite[Lemma 7.6]{Ev11} showed that the action
       $\Z=Symp_c(U)$ is free on its orbit in $\pi_0(\FF)$.  Therefore, the action of $\tau_L$ on $\pi_0(\FF)$ is
       also free and transitive, so actions of the generators $\tau_L$ and $\overline{l^{-1}}\circ\tau_\rho$ have to match,
      that is, $[\tau_L(C_+),\tau_L(C_-)]=[\overline{l^{-1}}\circ\tau_\rho(C_+), \overline{l^{-1}}\circ\tau_\rho(C_-)]\in\pi_0(\FF)$
      up to a change of the orientation of $L$. This shows that $\tau_L$ is Hamiltonian isotopic to $\overline{l^{-1}}\circ\tau_\rho$ in $Symp_c(U_2)$, since
      the stabilizer of the action of $Symp_c(U_2)$ on $\FF$ is weakly contractible by \eqref{e:ContF}.

\epf

\brmk\label{r:U2size} From the proof it is clear that the particular sizes of balls involved in Lemma \ref{l:BDlocal} is
not relevant, as long as the K\"ahler packing is possible.  By abuse of notation,
we will denote $U_2$ by open symplectic manifolds obtained this way.

\ermk

\brmk One may also construct a proof of Lemma \ref{l:BDlocal} from a more classical algebro-geometric point of view.  See
discussions in Section 5.
\ermk

We are now ready to return to the general case of $U$.  From \eqref{e:Br=Con} we have a $Br_n$-action on $\pi_0(\FF)$,
 which is free and transitive.  By comparing with the free action of $\pi_0(Symp_c(U_2))$ on the same space,
Evans obtains the monomorphism in Theorem \ref{t:EvIn}:
\beq\label{e:Defe} e:\pi_0(Symp_c(U))\hookrightarrow Br_n. \eeq

This is precisely the identification taken in Lemma \ref{l:BDlocal} in the case of $U_2$.
Therefore, one may interpret Lemma \ref{l:BDlocal}
as showing that $e$ is an isomorphism for $n=2$.  The following result shows that $e$ is an isomorphism for any $n$,
and that $\pi_0(Symp_c(U))$ is generated by Dehn twists along matching cycles, hence implying Proposition \ref{p:MainS}.

\blem\label{l:key} Let $T\subset \pi_0(Symp_c(U))$ be the subgroup generated by Lagrangian Dehn twists of matching cycles.
  Then there exist an isomorphism
      $\kappa:T\rightarrow Br_n$, such that the following diagram commutes:

\beq\label{e:Commutative} \xymatrixcolsep{5pc}\xymatrix{
T \ar@{->}[rd]^\sim_\kappa \ar@{^{(}->}[r] ^-c & \pi_0(Symp_c(U)) \ar@{->}[d]_e\\
&Br_n
 }
\eeq

\elem

To free up the notations we define:

\bdf\label{d:NI} Let $f: B(r)\subset\C^2\rightarrow (M,\w)$ be a symplectic embedding, and $\Sigma\subset M$ is a symplectic divisor.  Then
$f$ and $\Sigma$ is said to intersect \textit{normally} if $f^{-1}(\Sigma)=B(r)\bigcap \{z_2=0\}$.

\edf

\bpf[Proof of Lemma \ref{l:key}]
    We adopt Smith's model and notation from Section \ref{s:setup}.  Consider the blow-down of $M\backslash C_{n+1}$ along $C_i$, $i\leq n$.
    The blow-down is identified with a symplectic ball $B$ coming with the embedded balls $B_i$ resulted from $C_i$ for $i\leq n$.  Denote $\wt C_{n+2}\subset B$ as the proper transformation
     of $C_{n+2}$ under the blow-down.  As long as $r=\w(C_i)$ is sufficiently small,
    it is clear that one has an embedded symplectic ball $B_{i(i+1)}\hookrightarrow B$ for $1\leq i\leq n-1$, such that (see Figure 2):

    \bit
    \item $B_{i(i+1)}$ intersects $\widetilde C_{n+2}$ normally,
    \item $B_i\cup B_{i+1}\subset B_{i(i+1)}$,
    \item $B_{i(i+1)}\cap B_k=\emptyset$, for any $k\neq i,i+1$.
    \eit


For example, one may take $B_{i(i+1)}=B(\frac{2i-n}{(n+1)\sqrt{\pi}}, \frac{1}{(1+n)\sqrt{\pi}}+\frac{2}{\sqrt{N\pi}})$.
From the local construction Lemma \ref{l:BDlocal}, one sees that the action of each generator of the braid group
$\sigma_i$ on $\pi_0(\FF)$ is explicitly realized by a symplectomorphism coming from ball-swapping which is compactly supported in
$B_{i(i+1)}$.  Moreover, these symplectomorphisms are isotopic to Lagrangian Dehn twists of matching cycles supported in $B_{i(i+1)}$,
and such an identification is given precisely by $e\circ c$ in \eqref{e:Commutative} by the discussion
preceding Lemma \ref{l:key}.  Since both $e$ and $c$ are injective, all arrows in \eqref{e:Commutative} are indeed isomorphisms.

\epf

\section{A symplectomorphism by ball-swappings}\label{s:LagU}

\subsection{Constructing $\phi$ for Proposition \ref{p:MainL}}

In this section, we will prove Proposition \ref{p:MainL}.  The basic idea of constructing $\phi$ is again to use the ball-swapping
in the compactification.


Notice first that one may assume $L$ and $L'$ are homologous in $M=\CP^2\#n\overline{\CP^2}$.  This follows from
 the classification of homology classes of Lagrangian spheres in rational manifolds \cite[Theorem 1.4]{LW12}.  Since $L, L'$ are
disjoint from $C_{n+1}$ which is a line in $M$, their classes have to have the form of $[C_i]-[C_j]$ for $i,j\leq n$.
By certain Dehn twists along the standard Lagrangian spheres, one may assume
$[L]=[L']=[C_1]-[C_2]$ in $M$.

Recall also the following result from \cite{LW12}.

\bthm[\cite{LW12}, Theorem 1.1, 1.2]\label{t:MIT}
    Let $b^+(M)=1$ and $GT(A)\neq0$, $L$ be a Lagrangian sphere.  Then $A$ has an embedded representative with minimal intersection with $L$.  If $A$ is represented by an embedded sphere, then any given representatives of $A$ can be symplectically isotoped so that they achieve minimal intersections.
\ethm

Let $L,L'\subset U=M\backslash(C_{n+1}\cup C_{n+2})$, Theorem \ref{t:MIT} implies that one
may isotope $\{C_3,\dots,C_n\}$ to another standard configuration $\{S_3,\dots, S_n\}$, which are disjoint from $L$.
Extend this isotopy of spheres to a Hamiltonian isotopy
$\Psi_t$, then $\Psi_t^{-1}(L)$ isotopes $L$ away from $\{C_3,\dots, C_n\}$.  By performing the same type of isotopy to
$L'$, we may assume $\{C_3,\dots, C_n\}$ are disjoint from both $L$ and $L'$.

We now blow down along the set $\cE=\{C_3,\dots, C_n\}$ to obtain a set of balls $\cB=\{B_3,\dots, B_n\}$,
and the resulting manifold is denoted as $M_2=\CP^2\#2\overline{\CP^2}$, where we have the proper transformations of $L$ and
$L'$ denoted as $\widetilde{L}$ and $\widetilde{L}'$.  Note that by removing the proper transform
of $C_{n+1}$ and $C_{n+2}$ we have $M_2\backslash(\wt C_{n+1}\cup\wt C_{n+2})=U_2$.  One also have
from Remark \ref{r:U2size}:

\blem[\cite{BLW12}, Lemma A.1 and the discussions  following it]\label{t:unique}
    Lagrangian spheres in $U_2=M_2\backslash(\wt C_{n+2}\cup\wt C_{n+1})$ are unique up to Hamiltonian isotopy.
\elem

Combining discussions above,
there is a compactly supported Hamiltonian isotopy $\widetilde{\Phi}_t:U_2\rightarrow U_2$,
such that $\widetilde{\Phi}_t(\widetilde L)=\widetilde L'$
where $\widetilde{\Phi}_t=id$ near $\wt{C}_{n+1}\cup\wt C_{n+2}$.

\begin{figure}[h]
\centering
\includegraphics[trim = 0mm 20cm 0mm 0mm, clip, width=5.5in]{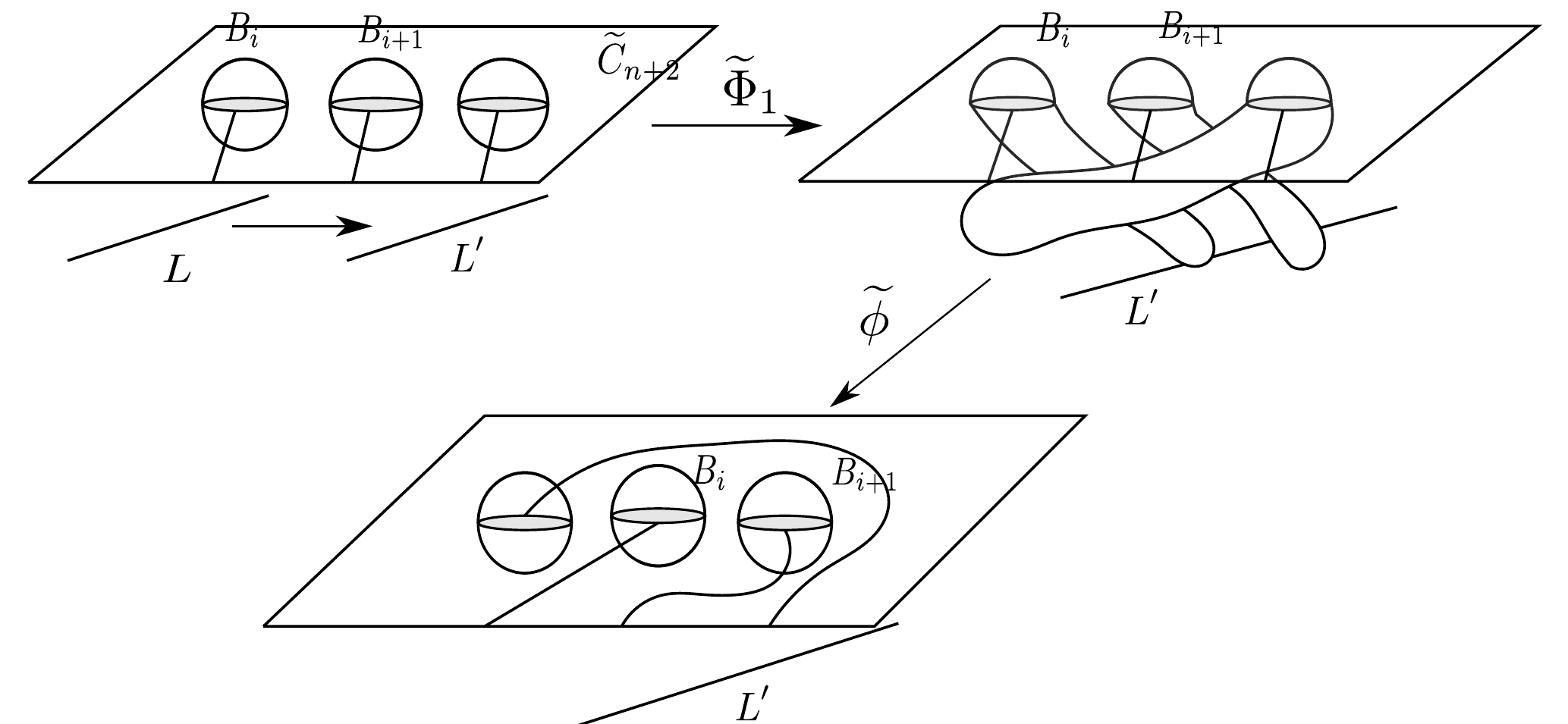}
\caption{A pictorial proof of Proposition \ref{p:MainL}} \label{f:pictorial proof}
\end{figure}

To obtain a compactly supported symplectomorphism on $M\backslash (C_{n+1}\cup C_{n+2})$,
it amounts to showing the following lemma, whose proof will be given in the next section.

\blem\label{l:RelConn}
   There is $\wt \phi\in Ham(M_2)$, such that:
   \ben[(i)]
   \item $\wt\phi(\wt \Phi_1(B_i))=B_i$ for $i\geq3$, $\wt \phi(\wt C_{n+2})=\wt C_{n+2}$,
   \item for some neighborhood $\NN$ of $\wt C_{n+1}\cup L'$,
         $\wt\phi|_{\NN}=id$.

   \een
\elem

Lemma \ref{l:RelConn} concludes the proof of Proposition \ref{p:MainL} because $\wt\phi\circ\wt\Phi$
sends $L$ to $L'$ and fixes $B_i$ for $i\geq3$, and it clearly is compactly supported in
$M_{2}\backslash \wt C_{n+1}$.  Therefore, it can be lifted to a compactly supported
symplectomorphism $\phi$ of $M\backslash C_{n+1}$ by blowing up the balls $B_i$ for $i\geq3$, which preserves $C_{n+2}$.
One then apply Lemma \ref{l:EvE} again to obtain a symplectomorphism $\phi$ which fixes $C_{n+2}$ pointwisely.  This
is always possible because $Ham_c(C_{n+2}\backslash C_{n+1})=Ham_c(D_2)\sim pt$ by the Smale's theorem.
Since the gauge group of the normal bundle of $C_{n+2}\backslash (C_{n+1}\cup C_{n+2})$ is homotopic equivalent
to $Map((S,*), SL_2(\R))\sim pt$, by composing another symplectomorphism fixing $C_{n+2}$ pointwisely one may
assume $\phi$ fixes also the normal bundle of $C_{n+2}$.  All these adjustments can be made
supported in a small neighborhood of $C_{n+2}$ thus not affecting $L'$.
Therefore, $\phi$ indeed descends to a compactly supported
symplectomorphism of $U=M\backslash (C_{n+1}\cup C_{n+2})$, which concludes our proof.

\brmk
For our purpose it suffices to prove $\wt\phi\in Symp(M_2)$ in Lemma \ref{l:RelConn}, which will be slightly easier.
From the proof it is clear that the lemma indeed holds for more general
cases of packing relative to a symplectic divisors, coupling with results in \cite{Mc12}
 and \cite{BLW12}, but we restrict ourselves for ease of expositions.
\ermk

\brmk
One recognizes that the braiding in the symplectomorphism group required for sending $L$ to $L'$
 comes exactly from the restriction of $\wt\phi$
on $C_{n+1}\backslash C_{n+2}$, which swaps the shadows $B_i\cap\wt C_{n+2}$ by $Diff(D_2)$.
 This nicely matches the pictures of Section \ref{s:MCGW}.
\ermk


\subsection{Connectedness of ball packing relative to a divisor}




We prove Lemma \ref{l:RelConn} in this section.  The overall idea is not new.  In the absence of $\wt C_{n+2}$, this is just
the ball-packing connectedness problem in the complement of a Lagrangian spheres in rational manifolds.  This was settled in
\cite{BLW12} using McDuff's non-generic inflation lemma \cite[Lemma 4.3.3]{Mc12} and her ideas in the original proof of
connectedness of ball-packings \cite[Section 3]{Mc98}.  We will adapt these ideas in our relative case and point out
necessary modifications.  We will continue to use notations defined in previous sections.


Let $\iota_0:B_i(r)\hookrightarrow M$ be the inclusion and $\iota_1=\wt \Phi_1:B_i\hookrightarrow M$.
 Fix some small $\delta>0$ and choose an extension of the embedding $\iota_j:B_i(r)\rightarrow M$ to $B_i(r+\delta)$ for
$j=0,1$, so that the extensions still intersect $\wt C_{n+2}$ normally.
Take two families of diffeomorphism $\phi_j^s$, $j=0,1$, so that the following holds:

\begin{enumerate}[(1)]
\item $\phi^0_j=id$,
\item $\phi_j^s|_{B(sr_i)}$ is a radial contraction from $B(sr_i)$ to $B(\delta)$, and identity near $\partial B(sr_i+\delta)$,
\item $\phi_j^s(\wt C_{n+2})=\wt C_{n+2}$.

\end{enumerate}

This is not hard to achieve if we have $B_i$ intersecting $\wt C_{n+2}$ normally in the first place.

The push-forward of $\omega$ by $\phi_j^s$ endows a family of symplectic forms
$\{\w_j^s\}_{0\leq s\leq 1}$ on $M$ for $j=0,1$.  Notice in our situation, $\wt\Phi_1$ in Lemma \ref{l:RelConn}
is compactly supported in $U_2$, which implies $\iota_0(B_i(\delta))=\iota_1(B_i(\delta))$ for
sufficiently small $\delta>0$. Therefore, performing a blow-up
on $\iota_j(B_i(\delta))$ gives a family of symplectic forms
$\{[\wt\w_j^s]\}$ on the same smooth manifold $M_2$ without further identification.
The resulting exceptional curves formed by blowing up
$B_i$ as $C_i$.  Notice $\int_{C_i}\wt\w^s_j=\pi(sr_i)^2$. We claim:

\blem\label{l:SITP}
    For any given $s\in [0,1]$, $\wt\w_0^s$ are isotopic to $\wt\w_1^s$ via a family of diffeomorphism
    $\psi^s(t):M\rightarrow M$ such that $\psi^s(t)(C_i)=C_i$ and $\psi^s(t)(L')=L'$ for
    any given $s\in[0,1]$ and $3\leq i\leq n+2$.
\elem

\bpf
    Notice that $\wt\w_0^s$ is deformation equivalent to $\wt\w_1^s$ by a family of symplectic forms
    $\{\Omega_t\}_{t\in[0,1]}$ preserving the orthogonality of $C_i$ with $C_{n+2}$ and constant near $L'$.
    A by-now-standard approach of correcting this
    deformation family into an isotopy of symplectic form, that is, a symplectic deformation by cohomologous symplectic forms, is done by inflation
    along $J$-holomorphic curves.




To put $L'$ into this framework, we may perform a symplectic cut near it, which
gives a symplectic $(-2)$ sphere $S$ (cf. \cite[Remark 2.2]{BLW12}).
The class of interests for inflation has the form
$[\tilde\w_0^s]-\epsilon\sum_{i=3}^n[C_i]$ with suitably chosen
$\epsilon$. When raised to a sufficiently large multiple, this class
has a nodal representative.  So far the argument is standard. To
inflate along this nodal curve, one considers $\cS$ as the union of
all its underlying irreducible components, and include $S$ and $C_i$
for $3\leq i\leq n+2$ if necessary.   Now we can follow the argument
in \cite[Proposition 1.2.9]{Mc12} and \cite[Lemma 1.1]{Mc13} for
this configuration $\cS$.  In principle, we need a family version of
inflation here, see discussion of \cite{Mc13} for details.  The
outcome is an isotopy of symplectic forms from $\wt\w_0^s$ to $\wt\w^0_1$ while keeping
$S$ and $C_i$ symplectic for $3\leq i\leq n+2$.


Lemma \ref{l:SITP} is then an immediate consequence of the relative Poincare
Lemma \cite[Theorem 2.3]{DJZ04} and the relative Moser technique \cite{MP94}.  In particular,
one first applies a family of diffeomorphism $g_t$ supported in a neighborhood of $\cup_{i=3}^{n+2}C_i$
preserving each of these divisors, so that $g_t^*(\wt \Omega_t)|_{C_{i}}=\wt \Omega_0|_{C_{i}}$.
Now the relative Poincare Lemma implies that $g_t^*(\wt \Omega_t)-\wt \Omega_0$ has a primitive
vanishing on $\cup_{i=3}^{n+2}C_i$ and a neighborhood of $L'$, where the Moser's technique implies the desired property in
Lemma \ref{l:SITP}.

\epf

To finish the proof of Lemma \ref{l:RelConn},
we already constructed a symplectomorphism $(\phi_1^s)^{-1}\psi^s\phi_0^s$
interpolating $\wt\w^s_0$ and $\wt\w^s_1$ for each $s$ from Lemma \ref{l:SITP}.
By blowing down $C_i$, one obtains a symplectomorphism $F^s$ sending
$\iota_0(B_i(sr_i))$ to $\iota_1(B_i(sr_i))$. From the smooth
dependence of the family $\psi^s$ on $s$, one learns that $\psi^1$
 is an isotopy of $M$.  By precomposing another isotopy $h$ preserving $C_i$ for all $3\leq i\leq n+1$
 as in the proof of \cite[Corollary 1.5]{Mc98}, one may further assume $\psi^1$ is
 identity near the $C_i$ for $i\leq n+1$.  One thus obtains a family of ball-packing by blowing
 down along $C_3,\dots,C_n$, which connects $\iota_0$ and $\iota_1$
 through a normal intersection family with $\wt C_{n+2}$ as desired.

\section{Concluding remarks on ball-swappings}

The ball-swapping technique we used in this paper seems to have rich
structures and could be of independent interests. This concluding
remark summarizes several possibly interesting directions of further
study of this class of objects.

\begin{itemize}
  \item It was pointed out to the author that,
        it seems instructive to compare ball-swappings with a closely
        related construction in algebraic geometry which is classical.  The author
        first learned about the following construction from Seidel's
        excellent lecture notes \cite{Se08}.  Consider $\C^n$ and
        its 2-point configuration space $\text{Conf}_2(\C^n)$.  One may
        associate to this space a fibration $\cE\rightarrow
        \text{Conf}_2(\C^n)$, where $\cE_b$ is the corresponding
        complex blow-up of at $b\in \text{Conf}_2(\C^n)$.  Seidel
        demonstrated in \cite[Example 1.12]{Se08} that, when $n=2$,
        one may partially compactify this family to $\overline\cE$
        by allowing the two points to collide, where the
        discriminant $\Delta$ is a smooth divisor.  A local normal
        disk $D$ centered at $p\in \Delta$ thus gives a
        sub-fibration over $D$, where the fiber over $p$ is a
        surface with a single ordinary double point.  Then \cite[Lemma
        1.11]{Se08} shows the monodromy around $\partial D$ is
        precisely a Dehn twist when an appropriate K\"ahler form is endowed on
        the fiber.  But this monodromy is equally clearly a
        ball-swapping, which establishes the relation between ball-swapping and
        the Dehn twists as monodromies in algebraic geometry.
        This algebro-geometric point of view also
        provides another (more elegant) proof for Lemma
        \ref{l:BDlocal}, except now one needs to sort out details for
        the last step of untwisting on the removed divisor.  In a
        Lefschetz fibration point of view, this untwist is
        equivalent to slowing down the Hamiltonian at infinity in the
        description of Seidel's Dehn twist.  One may consult
        \cite{Se97} for further details.

        However, the construction of ball-swapping is
        \textit{apriori} richer than Dehn twists even in dimension $4$.  Formally,
        given symplectic manifold $M^4$, consider the action of Hamiltonian group $Ham(M)$ on the space of
        ball-embeddings $Emb_\epsilon(M)=\{\phi:\coprod_{i=1}^{n}B(\epsilon_i)\rightarrow M\}$
        for $\epsilon=(\epsilon_1,\dots, \epsilon_n)$.  Take an
        orbit $\cO$ of the action, the stablizer of the action
        descends to a subgroup of $Symp(M\#n\overline{\CP^2})$.
        Then isotopy classes of ball-swappings are the images of
        $\pi_0(Symp(M\#n\overline{\CP^2}))$ under the connecting
        map from $\pi_1(Emb_\epsilon(M))$.  From this formal point
        of view, the above algebro-geometric construction
        amounts to the ball-swappings when restricted to images of

        $$\pi_1(\text{Conf}_n(M))\rightarrow\pi_1(Emb_\epsilon(M))\rightarrow \pi_0(Symp(M\#n\overline{\CP^2})),$$

        while the image of the second arrow forms the full ball-swapping subgroup.  In general, the inclusion
        $\text{Conf}_n(M)\rightarrow Emb_\epsilon(M)$
        is not a homotopy equivalence even when $n=1$ and $M$ is as simple
        as $S^2\times S^2$ with non-monotone forms \cite{LP04} (but is still $1$-connected!).  
        In this particular example one already sees the sizes of packed balls come
        into play.  In higher dimensions, the topology
        of space of ball-embedding in even $\C^n$ is completely
        open.   Therefore, it seems interesting to clarify the gap
        between the algebro-geometric construction and the full
        ball-swapping subgroup.

  \item The ball-swapping symplectomorphisms seem particularly useful in problems involving $\pi_0(Symp(M))$, when $M$ is
        a rational or ruled manifold.  From examples known to date and the algebro-geometric constructions
        above, it seems reasonable to speculate that ball-swappings in dimension $4$ is generally related to Lagrangian Dehn twists, at least for those with $b^+=1$.  A particular tempting question asks that, a blow-up at $n$ balls on $\CP^2$ with generic sizes has a connected
         symplectomorphism group (because they do not admit any Lagrangian Dehn twists, see \cite{LW12}).  But to reduce the subgroup generated by ball-swapping into problems of braids, as in what appeared in this note (and also \cite{Se08,Ev11}) seems to require independent efforts in finding an appropriate symplectic spheres, as well as a good control on the bubbles.

        Here we roughly sketch a viewpoint through ball-swapping to $\pi_0(Symp(M))$ when $M=\CP^2\#5\overline{\CP^2}$ where all extra technicalities can be proved irrelavant.   In this case, Evans \cite{Ev11} showed that $\pi_0(Symp(M))=Diff(S^2,*)$, where $*$ is a fixed set on $S^2$ consisting of $5$ points, which can be identified with a braid group on $S^2$.  This coincides with an observation of Seidel \cite[Example 1.13]{Se08}.
        Jonny Evans also explained to the author that these braids can be seen to be generated by Lagrangian Dehn twists, following from results of \cite{Ev11}\footnote{Private communications.}.

        From the ball-swapping point of view, by blowing down the 5 exceptional spheres of classes $E_1,\dots, E_5$, we have 5 embedded balls $B_1,\dots, B_5$ intersecting a $2H$-sphere $C$ normally.  Then the restriction of swapping these $5$ balls on the $2H$-sphere gives precisely a copy of the group $Diff(S^2,*)$, where the 5 points now is actually the 5 disks coming from intersections $B_i\cap C$.
        To get the actual statement of Evans, one still needs to switch between blow-ups and downs and go through Evans's proof to show all contributions from other components (configuration space of curves involved, automorphisms of normal bundles, etc.) cancels.

        The ball-swapping also gives a way of seeing these braidings actually come from Lagrangian Dehn
        twists, but we will not give full details.  A naive attempt is to include two embedded balls above to a larger ball normally intersecting $C$ as in the proof of Lemma \ref{l:key}.  However, this cannot work for packing size restrictions.  Instead, Evans shows that there is
        an embedded Lagrangian $\RP^2$ in the complement of a configuration of symplectic spheres consisting of classes $\{2H, E_1,\dots, E_5\}$.
        By removing this $\RP^2$, one can show that the remainder of $M$ can be identified with a symplectic fibration over $C$, where a generic
        fiber is a disk, and there are 5 singular fibers consisting of the union of a $E_i$-sphere and disk.
        One then choose a disk in $C$ separating the intersections of the $E_1$ and $E_2$-spheres with $C$ with other $E_i$'s, then the fibration restricted to this disk is a product of two $D_2$'s blown-up 2 points, denoted as $V_2$.  At this point one easily reproduces the proof of Lemma \ref{l:key} with an identification of the the symplectization of $V_2\backslash C$ with $T^*S^2$.

  \item In higher dimensions ball-swapping seems a new way of constructing monodromies and quite symplectic in flavor.  For example, when
        the ball-packing is sufficiently small, because of Darboux theorem one could essentially move the embedded balls as if they were just points. The question when such symplectomorphisms are actually Hamiltonian seems intriguing.

        One particular situation in question is when there is only one embedded ball.  Suppose a small embedded ball moves along a loop which is nontrivial in $\pi_1(M)$, is it true that the resulting ball-swapping always lies outside the Hamiltonian group?  \cite{Sm12} explained
        the existence of a Lagrangian torus near a small blow-up.  It is not difficult to find a symplectomorphism in the component of $Symp(M)$
        where the ball-swapping lies, so that this torus is invariant under it.  Thus, in the sense of \cite{Se11}, one constructs a family of objects in the Fukaya category of the mapping torus of $M\#\overline{\CP^2}$ naturally associated to a ball-swapping.

  \item In principle, one may extend the ``swapping philosophy" to a more general background of surgeries to obtain automorphisms, even in non-symplectic situations.  As a simple example, consider a Hamiltonian loop of Lagrangian $2$-torus in $M^4$, by which we mean a path of Hamiltonians $\phi_t\subset Ham(M)$ such
      that $\phi_{0,1}(L)=L$.  One may perform a Luttinger surgery on $L$ \cite{Lu95}, which gives a new symplectic manifold $M^L$.  Then the ``swapping" $\phi_1$ is lifted to an automorphism of $M^L$.
         
         One may also choose any symplectic neighborhood of a symplectic submanifold, as long as one has interesting loops of these objects.  Naively thinking, when the submanifold is actually a divisor, the swapping symplectomorphism should always be Hamiltonian, but the author has no proof for that.  Otherwise, since blowing up a divisor only gives a deformation of the symplectic form, one may always deform the form first and then swap the divisor.  This will yield a symplectomorphism of the original manifold $M$, which should contain interesting information if not a Hamiltonian.

\end{itemize}




\begin{thebibliography}{99}



\bibitem{Ab12}
M.~Abouzaid.
\newblock Nearby Lagrangians with vanishing Maslov class are homotopy equivalent.
\newblock {\em Invent. Math.} 189 (2012), no. 2, 251-313.

\bibitem{AI12}
M.~Abouzaid, I.~Smith.
\newblock Exact Lagrangians in plumbings.
\newblock preprint http://arxiv.org/abs/1107.0129





\bibitem{BLW12}
M.S.~Borman, T.-J.~Li, W.~Wu.
\newblock Spherical Lagrangians via ball packings and symplectic cutting
\newblock {\em To appear in Selecta Math.}, http://arxiv.org/abs/1211.5952



\bibitem{DJZ04}
W.~Domitrz, W.~Janeczko, M.~Zhitomirskii.
\newblock Relative Poincare lemma, contractibility, quasi-homogeneity and vector fields tangent to a singular variety.
\newblock {\em Illinois J. Math.} 48 (2004), no. 3, 803-835.



\bibitem{EvT}
J.D.~Evans.
\newblock Symplectic topology of some Stein and rational surfaces
\newblock {\em Ph. D. Thesis}, 2010 Christ's College, Cambridge.

\bibitem{Ev10}
J.D. Evans.
\newblock Lagrangian spheres in del {P}ezzo surfaces.
\newblock {\em J. Topol.}, 3(1):181--227, 2010.

\bibitem{Ev11}
J.D.~Evans.
 \newblock Symplectic mapping class groups of some Stein and rational surfaces.
 \newblock {\em J. Symplectic Geom.}, 9 (2011), no. 1, 45-82.



\bibitem{Go95}
R.E.~ Gompf.
\newblock A new construction of symplectic manifolds.
\newblock {\em Ann. of Math. (2)}, 142(3):527--595, 1995.


\bibitem{Hi12}
R.~Hind.
\newblock Lagrangian unknottedness in Stein surfaces.
\newblock {\em Asian J. Math.} 16 (2012), no. 1, 1-36.

\bibitem{HPW12}
R.~Hind, M.~Pinsonnault, W.~Wu.
\newblock Symplectomorphism groups of non-compact manifolds and space of Lagrangians.
\newblock preprint.

\bibitem{IU05}
A.~Ishii, H.~Uehara.
\newblock Autoequivalences of derived categories on the minimal resolutions of $A_n$ singularities on surfaces.
\newblock {\em J. Differential Geom.} 71 (2005), 385-435

\bibitem{IUU10}
A. Ishii, K. Ueda, and H. Uehara.
\newblock Stability conditions on $A_n$ singularities.
\newblock {\em J. Differential Geom.} 84 (2010), no. 1, 87-126

\bibitem{KS02}
M.~Khovanov and P.~Seidel.
\newblock Quivers, Floer cohomology, and braid group actions.
\newblock {\em J. Amer. Math. Soc.} 15 (2002), no. 1, 203-271.

\bibitem{Lu95}
K. Luttinger
\newblock Lagrangian tori in $R^4$.
\newblock {\em J. Differential Geom.} 42 (1995), no. 2, 220¨C228. 

\bibitem{Sm12}
I.~Smith.
\newblock Floer cohomology and pencils of quadrics.
\newblock {\em Invent. Math.} 189 (2012), no. 1, 149-250.

\bibitem{Le95}
E.~Lerman.
\newblock Symplectic cuts.
\newblock {\em Math. Res. Lett.}, 2(3):247-258, 1995.



\bibitem{LL01}
T.-J.~Li, A.~Liu.
\newblock Uniqueness of symplectic canonical class, surface cone
and symplectic cone of 4-manifolds with $B^+=1$.
\newblock {\em J. Differential Geom.}, 58(2001), no. 2, 331-370.



\bibitem{LM12}
Y.~Lekili, M.~Maydanskiy.
\newblock The symplectic topology of some rational homology balls.
\newblock To appear in {\em Comm. Math. Helvetici.}, arxiv:1202.5625.

\bibitem{LP04}
F. ~Lalonde, M. ~Pinsonnault.
\newblock The topology of the space of symplectic balls in rational 4-manifolds.
\newblock {\em Duke Math. J.} 122 (2004), no. 2, 347¨C397.


\bibitem{LW12}
T.-J.~Li, W.~Wu.
\newblock Lagrangian spheres, symplectic surfaces and the
symplectic mapping class group.
\newblock {\em Geom. Topol.}, 16(2):1121-1169, 2012.

\bibitem{Mc98}
D.~McDuff.
\newblock From symplectic deformation to isotopy.
\newblock In {\em Topics in symplectic {$4$}-manifolds
({I}rvine, {CA}, 1996)},
First Int. Press Lect. Ser., I, pages 85-99. Int. Press,
Cambridge, MA, 1998.



\bibitem{Mc12}
D.~McDuff.
\newblock Nongeneric {$J$}-holomorphic curves in rational manifolds.
\newblock Preprint, 2012.

\bibitem{Mc13}
D.~McDuff.
\newblock Symplectic embeddings of 4-dimensional ellipsoids: erratum.
\newblock http://arxiv.org/abs/1305.0230.

\bibitem{MP94}
D.~McDuff, L.~Polterovich.
\newblock Symplectic packings and algebraic geometry.
\newblock {\em Invent. Math.}, 115(3):405--434, 1994.
\newblock With appendix by Y. Karshon.

\bibitem{MS98}
D.~McDuff, D.~Salamon.
\newblock Introduction to symplectic topology. Second edition.
\newblock Oxford Mathematical Monographs. The Clarendon Press, Oxford University Press, New York, 1998. x+486 pp. ISBN: 0-19-850451-9



\bibitem{Pi08}
M.~Pinsonnault.
\newblock Maximal compact tori in the Hamiltonian group of 4-dimensional symplectic manifolds.
\newblock {\em J. Mod. Dyn.} 2 (2008), no. 3, 431-455.

\bibitem{Ri10}
A.~Ritter.
\newblock Deformations of Symplectic Cohomology and Exact Lagrangians in ALE Spaces.
\newblock {\em Geom. Funct. Anal.} 20 (2010), no. 3, 779-816.

\bibitem{Se97}
P.~Seidel.
\newblock Floer homology and the symplectic isotopy problem.
\newblock Thesis, {\em University of Oxford, 1997}.

\bibitem{Se98}
P.~Seidel.
\newblock Symplectic automorphisms of $T^*S^2$
\newblock {\em http://arxiv.org/abs/math/9803084.}


\bibitem{Se99}
P.~Seidel.
\newblock Lagrangian two-spheres can be symplectically knotted.
\newblock {\em J. Differential Geom.} 52 (1999), no. 1, 145-171.

\bibitem{Se08}
P.~Seidel.
\newblock Lectures on four-dimensional {D}ehn twists.
\newblock In {\em Symplectic 4-manifolds and algebraic
surfaces}, volume 1938
of {\em Lecture Notes in Math.}, pages 231-267. Springer,
Berlin, 2008.

\bibitem{Se11}
P.~Seidel.
\newblock Abstract analogues of flux as symplectic invariants.
\newblock http://arxiv.org/abs/1108.0394



\bibitem{Se12}
P.~Seidel.
\newblock Lagrangian homology spheres in $(A_m)$ Milnor fibres via $C^*$-equivariant $A_\infty$ modules.
\newblock {\em Geom. Topo.} 16 (2012), 2343-2389



\end{thebibliography}

\end{document}